\documentclass[12pt]{article}
\usepackage{amsmath,amsfonts,latexsym,amsthm,amssymb}
\newcommand{\D}{\mathbb D_t^{(\alpha )}}
\newcommand{\Rn}{\mathbb R^n}
\newcommand{\CN}{\mathbb C^N}
\numberwithin{equation}{section}

\DeclareMathOperator{\const}{const}

\DeclareMathOperator{\R}{Re}
\begin{document}
\newtheorem{prop}{Proposition}
\newtheorem{lem}{Lemma}
\newtheorem{teo}{Theorem}
\pagestyle{plain}
\title{Fractional-Parabolic Systems}
\author{Anatoly N. Kochubei\footnote{This work was supported in part by the Ukranian Foundation for Fundamental Research under
Grant 28.1/017.}
\\ \footnotesize Institute of Mathematics,\\
\footnotesize National Academy of Sciences of Ukraine,\\
\footnotesize Tereshchenkivska 3, Kiev, 01601 Ukraine\\
\footnotesize E-mail: \ kochubei@i.com.ua}
\date{}
\maketitle

\vspace*{3cm}
\begin{abstract}
We develop a theory of the Cauchy problem for linear evolution systems of partial differential equations with the Caputo-Dzhrbashyan fractional derivative in the time variable $t$. The class of systems considered in the paper is a fractional extension of the class of systems of the first order in $t$ satisfying the uniform strong parabolicity condition. We construct and investigate the Green matrix of the Cauchy problem. While similar results for the fractional diffusion equations were based on the H-function representation of the Green matrix for equations with constant coefficients (not available in the general situation), here we use, as a basic tool, the subordination identity for a model homogeneous system. We also prove a uniqueness result based on the reduction to an operator-differential equation.
\end{abstract}
\vspace{2cm}
\newpage
\section{INTRODUCTION}

Fractional diffusion equations of the form
\begin{equation}
\left( \D u\right) (t,x)-Au(t,x)=f(t,x),\quad 0\le t\le T,x\in \Rn,
\end{equation}
where $0<\alpha <1$, $\D$ is the Caputo-Dzhrbashyan fractional derivative, that is
$$
\left( \D u\right) (t,x)=\frac{1}{\Gamma (1-\alpha )}\left[
\frac{\partial}{\partial t}\int\limits_0^t(t-\tau )^{-\alpha
}u(\tau ,x)\,d\tau -t^{-\alpha }u(0,x)\right] ,
$$
$A$ is a second order elliptic operator, are among the basic subjects in the theory of fractional differential equations. The initial motivation came from physics -- the equations of the above type were first used for modeling anomalous diffusion on fractals by Nigmatullin \cite{Nigm} and for a description of Hamiltonian chaos by Zaslavsky \cite{Zasl}. See the survey papers \cite{GM,MK} for a description of the present status of this research area.

The first mathematical works in this direction dealt either with the case of an abstract operator $A$, that is with a kind of an abstract Cauchy problem \cite{K1} (see \cite {KST} for further references), or with the case where $A=\Delta$ is a Laplacian. For the latter case, a fundamental solution of the Cauchy problem (FSCP) is expressed via Fox's H-function \cite{K2,SW}; uniqueness theorems were proved in \cite{K2} for an equation with a general second order elliptic operator $A$; see also \cite{Ps}. The first example of an initial-boundary value problem for this equation was considered in \cite{W}. Later the initial-boundary value problems for fractional diffusion equations were studied in \cite{MNV}, with an emphasis on the probabilistic aspects, and in \cite{L}; for the probabilistic interpretations see also \cite{S,MS} and references therein.

In \cite{EK} (see also \cite{EIK}), Eidelman and the author constructed and investigated a FSCP for fractional diffusion equations with variable coefficients. We followed the classical parametrix method using an H-function representation for the parametrix kernel and the detailed information about asymptotic properties of the H-function available from \cite{B,KS}. The fractional diffusion equation shares many essential properties with second order para\-bolic equations (though some properties are different, like, for example, the singularity of the FSCP at $x=0$ appearing for $n\ge 2$).

In the development of the theory of linear partial differential equations of parabolic type, the next step after the study of second order equations was to identify a class of systems which can be called parabolic. For such systems, there must be a well-posed Cauchy problem whose fundamental solution is an ordinary function smooth outside the singular point $t=0$. Such a class of systems was first found by Petrowsky \cite{P} in 1938. In fact, after Petrowsky's work it was understood that a complete theory of partial differential equations must include a thorough study of systems of equations. Petrowsky introduced and investigated not only parabolic systems but also hyperbolic ones, systems with correct Cauchy problems etc. For several decades, these subjects were central in the theory of partial differential equations; see, in particular, \cite{GS3,Pal}. The class of parabolic systems was studied in the greatest detail  \cite{E,EIK,F,LSU}.

So far only some special cases of partial fractional-differential systems have been considered \cite{GD,GDM,Pietr,VV}. However it is the author's opinion that the fractional calculus is mature enough to initiate a general theory of systems of partial fractional-differential equations. In this paper we follow the above line and find a fractional analog of the class of parabolic systems. Note that simple examples of such fractional parabolic systems appear as linearized two-component fractional reaction-diffusion systems used in the study of self-organization phenomena; see \cite{GD,GDM} and references therein.

We consider systems of the form (1.1) where $u=(u_1,\ldots ,u_N)$ is a vector-valued function,
\begin{equation}
A=A(x,D_x)=A_0(x,D_x)+A_1(x,D_x),
\end{equation}
is a differential operator of even order $2b$ with matrix-valued coefficients,
\begin{equation}
(A_0(x,D_x)u)_i=\sum\limits_{j=1}^N\sum\limits_{|\beta |=2b}a_\beta ^{ij}(x)D_x^\beta u_j,\quad i=1,\ldots ,N,
\end{equation}
\begin{equation}
(A_1(x,D_x)u)_i=\sum\limits_{j=1}^N\sum\limits_{|\beta |<2b}a_\beta ^{ij}(x)D_x^\beta u_j,\quad i=1,\ldots ,N,
\end{equation}
$$
D_x^\beta =D_{x_1}^{\beta_1}\cdots D_{x_n}^{\beta_n},\quad D_{x_j}=\frac1{\sqrt{-1}}\frac{\partial }{\partial x_j},\quad |\beta |=\beta_1+\cdots +\beta_n.
$$
We assume that all the coefficients $a_\beta ^{ij}(x)$ are bounded and satisfy the global H\"older condition
$$
\left| a_\beta ^{ij}(x)-a_\beta ^{ij}(y)\right| \le C|x-y|^\gamma
$$
(below the letters $C,c$ will denote various positive constants while $\gamma >0$ will denote all the H\"older exponents; for simplicity, we denote by $|\cdot |$ norms of all finite vectors and matrices). We also assume the uniform strong parabolicity condition: for all $\eta \in \Rn$, $z\in \mathbb C^N$,
\begin{equation}
\R \langle A_0(x,\eta )z,z\rangle \le -\delta |\eta |^{2b}|z|^2,\quad \delta >0.
\end{equation}

In fact, the key ingredients in the construction of a FSCP for a problem with variable coefficients are precise estimates for the model problem
\begin{equation}
\left( \D u(t,x)\right) =A_0(y,D_x)u(t,x)
\end{equation}
containing only the homogeneous highest order differential operator in $x$, with ``frozen'' coefficients depending on a parameter point $y$. As the first step, one has to consider the case where the coefficients $a_\beta ^{ij}(x)$, $|\beta |=2b$, are constant. Already in this case, the study of a FSCP is far from trivial. The approach used in \cite{EK,EIK} based on the H-function representation, does not work for systems.

Instead, we use the subordination representation \cite{Ba,BM} expressing the FSCP for the model system via the FSCP for the first-order (in $t$) parabolic system. At the first sight, it looks an easy approach to all fractional problems. However the subordination identity involves the integration in $t$ over the half-axis $(0,\infty )$, while usually a FSCP for a parabolic equation or system is constructed only on a finite time interval. Nevertheless, for our model case of constant coefficients the subordination method works efficiently giving, by the way, new proofs of the estimates known for fractional diffusion equations. Note also that the probabilistic side of subordination, not touched here, is also an important subject of fractional analysis; see \cite{BMN,Kol,MTM,OB}.

We also prove a uniqueness theorem for general systems (1.1). Again, the method of proving uniqueness in \cite{K2,EIK} (based on a kind of the maximum principle) is applicable only for second order equations. Here we use the reduction to an abstract equation from \cite{K1} and the regularized resolvent estimate from \cite{HHN}.

The main results of this paper are collected in Section 2. Section 3 contains miscellaneous auxiliary results used subsequently. Proofs of the estimates for the estimates for the Green matrix of a homogeneous system with constant coefficients are given in Section 4 and are complemented in Section 5 with some considerations regarding the parametrix kernels. In Section 6, these results are used to substantiate the Levi method in our situation. The proof of the uniqueness theorem is given in Section 7.

The author is grateful to the anonymous referee for helpful comments and suggestions.

\medskip
\section{Main results}

In this section we introduce basic notions and formulate principal results. The proofs will be given in subsequent sections.

{\it 2.1. The model system. Subordination}. Let us consider systems of the form
\begin{equation}
\left( \D u(t,x)\right) =A_0(D_x)u(t,x)
\end{equation}
where
$$
(A_0(D_x)u)_i=\sum\limits_{j=1}^N\sum\limits_{|\mu |=2b}a_\mu ^{ij}D_x^\mu u_j,\quad i=1,\ldots ,N,
$$
$a_\mu ^{ij}\in \mathbb C$, and for any $\eta \in \Rn$, $z\in \mathbb C^N$,
\begin{equation}
\R \langle A_0(\eta )z,z\rangle \le -\delta |\eta |^{2b}|z|^2,\quad \delta >0.
\end{equation}

Under the assumption (2.2) (in fact, even under a weaker assumption of parabolicity in the sense of Petrowsky), the differential expression $A_0(D)$ defines on the space $L^2(\Rn ,\mathbb C^N)$ of square integrable vector-functions with values in $\mathbb C^N$, an infinitesimal generator $\mathcal A_0$ of a $C_0$-semigroup $S_1(t)=e^{t\mathcal A_0}$ (see \cite{Kr}).

By the subordination theorem (see Theorem 3.1 in \cite{Ba}), the system (2.1) interpreted as an equation in $L^2(\Rn ,\mathbb C^N)$, possesses a solution operator $S_\alpha (t)$, such that for any element $u_0 =u_0 (x)$ from the domain $D(\mathcal A_0)$, the function $u(t,x)=(S_\alpha (t)u_0 )(x)$, $t\ge 0$, $x\in \Rn$, is a solution of the equation (2.1) satisfying the initial condition $u(0,x)=u_0 (x)$. In addition,
\begin{equation}
S_\alpha (t)=\int\limits_0^\infty \varphi_{t,\alpha }(s)S_1(s)\,ds,\quad t\ge 0,
\end{equation}
where $\varphi_{t,\alpha }(s)=t^{-\alpha }\Phi_\alpha (st^{-\alpha})$,
$$
\Phi_\alpha (\zeta )=\sum\limits_{k=0}^\infty \frac{(-\zeta )^k}{k!\Gamma (-\alpha k+1-\alpha )},
$$
so that $\Phi_\alpha$ can be written as the Wright function
\begin{equation}
\Phi_\alpha (\zeta )={}_0\Psi_1\Bigl[ \begin{matrix}
- \\ (1-\alpha ,-\alpha )\end{matrix}\Bigl| -\zeta \Bigr] .
\end{equation}
See \cite{KST} for general information regarding the definition and properties of the Wright functions; see also Section 3.3 below.

The function $\Phi_\alpha$ is a probability density:
$$
\Phi_\alpha (t)\ge 0,\ t>0;\quad \int\limits_0^\infty \Phi_\alpha (t)\,dt=1.
$$
It is connected also with the Mittag-Leffler function
\begin{equation}
E_\alpha (\zeta )=\sum\limits_{k=0}^\infty \frac{\zeta^k}{\Gamma (1+\alpha k)},\quad \zeta \in \mathbb C,
\end{equation}
via the Laplace transform identity
\begin{equation}
E_\alpha (-\zeta )=\int\limits_0^\infty \Phi_\alpha (t)e^{-\zeta t}\,dt,\quad \zeta \in \mathbb C.
\end{equation}

By the classical theory of parabolic equations, the semigroup $S_1(t)$ possesses the integral representation
$$
(S_1(t)\eta )(x)=\int\limits_{\Rn}Z(t,x-\xi )\eta (\xi )\,d\xi ,\quad \eta \in L^2(\Rn ,\mathbb C^N),
$$
in terms of the FSCP $Z(t,x)$ of the parabolic system $\dfrac{\partial u}{\partial t}=A_0(D_x)u$. It follows from the estimates of $\Phi_\alpha$ and $Z$ (see below) that, for example, if $\eta \in \mathcal S(\Rn )$, then
$$
(S_\alpha (t)\eta )(x)=\int\limits_{\Rn}Z_\alpha (t,x-\xi )\eta (\xi )\,d\xi
$$
where
\begin{equation}
Z_\alpha (t,x)=\int\limits_0^\infty \varphi_{t,\alpha }(s)Z(s,x)\,ds,\quad x\ne 0
\end{equation}
(as we have seen for the diffusion equations \cite{EK,EIK}, $Z_\alpha$ may have a singularity at $x=0$).

The kernel $Z_\alpha$ is a FSCP for the system (2.1).

In order to obtain an integral representation of a solution $u(t,x)$ of the inhomogeneous equation
$$
\left( \D u(t,x)\right) -A_0(D_x)u(t,x)=f(t,x),\quad u(0,x)=u_0(x),
$$
in the form
$$
u(t,x)=\int\limits_{\Rn}Z_\alpha (t,x-\xi )u_0(\xi )\,d\xi +\int\limits_0^td\tau \int\limits_{\Rn}Y_\alpha (t-\tau ,x-y)f(\tau ,y)\,dy
$$
(the definition of a classical solution will be given below for a more general situation), we need another kernel
$$
Y_\alpha (t,x)=\left( \mathbb D_t^{(1-\alpha )}Z_\alpha \right) (t,x),\quad x\ne 0.
$$
As we will see,
\begin{equation}
Y_\alpha (t,x)=\int\limits_0^\infty \psi_{t,\alpha }(s)Z(s,x)\,ds,
\quad x\ne 0
\end{equation}
where
\begin{equation}
\psi_{t,\alpha }(s)=t^{-1}{}_0\Psi_1\Bigl[ \begin{matrix}
- \\ (0,-\alpha )\end{matrix}\Bigl| -st^{-\alpha} \Bigr] .
\end{equation}

Note that $\varphi_{t,\alpha }$ is the law of the inverse to a $\alpha$-stable subordinator; see \cite{MS}. In the context of fractional calculus models of physical phenomena, the functions $\varphi_{t,\alpha }$ and $\psi_{t,\alpha }$ appeared for the first time in \cite{Main}.

Using properties of the functions $Z$ and $\Phi_\alpha$ we get the integral identities
$$
\int\limits_{\Rn}Z_\alpha (t,x)\,dx=1,\quad \int\limits_{\Rn}Y_\alpha (t,x)\,dx=\frac{t^{\alpha -1}}{\Gamma (\alpha)}.
$$

\medskip
\begin{teo}
The matrix-functions $Z_\alpha (t,x)$, $Y_\alpha (t,x)$ are infinitely differentiable for $t>0$, $x\ne 0$, and satisfy the following estimates. Denote
$$
R=t^{-\alpha }|x|^{2b},\quad \rho (t,x)=\left( t^{-\alpha }|x|^{2b}\right)^{\frac1{2b-\alpha}}.
$$
\begin{description}
\item[(i)] If $R\ge 1$, then
\begin{equation}
\left| D_x^\beta Z_\alpha (t,x)\right| \le Ct^{-\alpha \frac{ (n+|\beta |)}{2b}}\exp (-\sigma \rho (t,x)),\quad \sigma >0;
\end{equation}
\begin{equation}
\left| D_x^\beta Y_\alpha (t,x)\right| \le Ct^{-\alpha \frac{ (n+|\beta |)}{2b}+\alpha -1}\exp (-\sigma \rho (t,x)).
\end{equation}
\item[(ii)] If $R\le 1$, $n+|\beta |<2b$, then
\begin{equation}
\left| D_x^\beta Z_\alpha (t,x)\right| \le Ct^{-\alpha\frac{ (n+|\beta |)}{2b}},
\end{equation}
\begin{equation}
\left| D_x^\beta Y_\alpha (t,x)\right| \le Ct^{-\alpha\frac{ (n+|\beta |)}{2b}+\alpha -1}.
\end{equation}
\item[(iii)] If $R\le 1$, $n+|\beta |>2b$, then
\begin{equation}
\left| D_x^\beta Z_\alpha (t,x)\right| \le Ct^{-\alpha}|x|^{-n+2b-|\beta |},
\end{equation}
\begin{equation}
\left| D_x^\beta Y_\alpha (t,x)\right| \le Ct^{-1}|x|^{-n+2b-|\beta |}.
\end{equation}
\item[(iv)] If $R\le 1$, $n+|\beta |=2b$, then
\begin{equation}
\left| D_x^\beta Z_\alpha (t,x)\right| \le Ct^{-\alpha},\quad \text{if $n=1$};
\end{equation}
\begin{equation}
\left| D_x^\beta Z_\alpha (t,x)\right| \le Ct^{-\alpha}[|\log (t^{-\alpha }|x|^{2b})|+1],\quad \text{if $n\ge 2$};
\end{equation}

\begin{equation}
\left| D_x^\beta Y_\alpha (t,x)\right| \le Ct^{-1}.
\end{equation}
\item[(v)] If $R\ge 1$, then
\begin{equation}
\left| \frac{\partial Z_\alpha (t,x)}{\partial t} \right|\le Ct^{- \frac{ \alpha n}{2b}-1}\exp (-\sigma \rho (t,x)).
\end{equation}
\item[(vi)] If $R\le 1$, $n<2b$, then
\begin{equation}
\left| \frac{\partial Z_\alpha (t,x)}{\partial t} \right|\le Ct^{- \frac{ \alpha n}{2b}-1}.
\end{equation}

If $R\le 1$, $n>2b$, then
\begin{equation}
\left| \frac{\partial Z_\alpha (t,x)}{\partial t} \right|\le Ct^{- \alpha -1}|x|^{-n+2b}.
\end{equation}

If $R\le 1$, $n=2b$, then
\begin{equation}
\left| \frac{\partial Z_\alpha (t,x)}{\partial t} \right|\le Ct^{- \alpha -1}[|\log (t^{-\alpha }|x|^{2b})|+1].
\end{equation}
\end{description}

In all the above estimates, the constants depend only on $N,n$, $\max \left| a_\mu^{ij}\right|$, and the strong parabolicity constant $\delta$.
\end{teo}

\bigskip
The estimates (2.10)-(2.22) agree with their counterparts for the fractional diffusion equations \cite{EK,EIK}, though in the latter case, for some values of $n$, there are more precise estimates of $Z_\alpha$ and $Y_\alpha$.

Note also that the fractional derivative $\D Z_\alpha$ satisfies the same estimate as the derivatives $D_x^\beta Z_\alpha$, $|\beta |=2b$.

\medskip
{\it 2.2. The general case}. As stated in Introduction, we consider the system (1.1)--(1.4) with bounded H\"older continuous coefficients, under the uniform strong parabolicity condition (1.5).

We call a vector-function $u(t,x)$, $0\le t\le T$, $x\in \Rn$, {\it a classical solution} of the system (1.1), with the initial condition
\begin{equation}
u(0,x)=u_0(x),\quad x\in \Rn ,
\end{equation}
if:
\begin{description}
\item[(i)]
$u(t,x)$ is continuously differentiable in $x$ up to the order $2b$, for each $t>0$;
\item[(ii)]
for each $x\in \Rn$, $u(t,x)$ is continuous in $t$ on $[0,T]$, and its fractional integral
\begin{equation}
\left( I_{0+}^{1-\alpha }u\right) (t,x)=\frac{1}{\Gamma (1-\alpha
)}\int\limits_0^t(t-\tau )^{-\alpha }u(\tau ,x)\,d\tau
\end{equation}
is continuously differentiable in $t$ for $0\le t\le T$.
\item[(iii)]
$u(t,x)$ satisfies the equation (1.1) and the initial condition (2.23).
\end{description}

\medskip
A classical solution $u(t,x)$ is called {\it a uniform classical solution}, if it is continuous in $t$ uniformly with respect to $x\in \Rn$, and the first derivative of the fractional integral (2.24) exists uniformly with respect to $x\in \Rn$.

Our main task is to construct {\it a Green matrix} for the problem (1.1), (2.23), that is such a pair
$$
\left\{ Z_\alpha^{(1)}(t,x;\xi ),Y_\alpha^{(1)}(t,x;\xi )\right\} ,\quad t\in [0,T],\ x,\xi \in \Rn ,
$$
that for any bounded function $f$, jointly continuous in $(t,x)$ and locally H\"older continuous in $x$ uniformly with respect to $t$, and any bounded locally H\"older continuous function $u_0$, the function
\begin{equation}
u(t,x)=\int\limits_{\mathbb R^n} Z_\alpha^{(1)}(t,x;\xi )u_0(\xi )\,d\xi+\int\limits_0^td\lambda \int\limits_{\mathbb R^n} Y_\alpha^{(1)}(t-\lambda,x;y)f(\lambda ,y)\,dy.
\end{equation}
is a classical solution of the problem (1.1),(2.23).

Denote by $Z_\alpha^{(0)}(t,x-\xi ;y)$ and $Y_\alpha^{(0)}(t,x-\xi ;y)$ the kernels defined just as  $Z_\alpha (t,x-\xi )$ and $Y_\alpha (t,x-\xi )$, but for the system (1.6) with the coefficients $a_\beta^{ij}$, $|\beta |=2b$, ``frozen'' at a point $y\in \Rn$, and other coefficients set equal to zero; in (1.6), $y$ appears as a parameter.

\medskip
\begin{teo}
(a) There exists a Green matrix  $\left\{ Z_\alpha^{(1)}(t,x;\xi ),Y_\alpha^{(1)}(t,x;\xi )\right\}$ of the form
$$
Z_\alpha^{(1)}(t,x;\xi )=Z_\alpha^{(0)}(t,x-\xi ;\xi )+V_Z(t,x;\xi ),
$$
$$
Y_\alpha^{(1)}(t,x;\xi )=Y_\alpha^{(0)}(t,x-\xi ;\xi )+V_Y(t,x;\xi ),
$$
where the kernels $Z_\alpha^{(0)}(t,x;\xi )$, $Y_\alpha^{(0)}(t,x;\xi )$ satisfy the estimates listed in Theorem 1 with coefficients independent on the parameter point $\xi$. The functions $V_Z$, $V_Y$ satisfy the following estimates.

\begin{description}
\item[(i)]
If $n+|\beta |<2b$, then
\begin{equation}
\left| D_x^\beta V_Z(t,x;\xi )\right| \le Ct^{-\frac{\alpha}{2b}(|\beta |+\gamma_0)}|x-\xi |^{-n+\gamma -\gamma_0}e^{-\sigma \rho (t,x-\xi )},\quad 0<\gamma_0<\gamma ,\ \sigma >0;
\end{equation}
\begin{equation}
\left| D_x^\beta V_Y(t,x;\xi )\right| \le Ct^{-1+\alpha -\frac{\alpha |\beta|}{2b}}|x-\xi |^{-n+\gamma }e^{-\sigma \rho (t,x-\xi )}.
\end{equation}

\item[(ii)]
If $n+|\beta |\ge 2b$, $|\beta |<2b,$ then
\begin{equation}
\left| D_x^\beta V_Z(t,x;\xi )\right| \le Ct^{-\alpha +\frac{\alpha \gamma_0}{2b}}|x-\xi |^{-n+2b-|\beta |+\gamma -\gamma_0}e^{-\sigma \rho (t,x-\xi )};
\end{equation}
\begin{equation}
\left| D_x^\beta V_Y(t,x;\xi )\right| \le Ct^{-1+\frac{\alpha \gamma_0}{b}}|x-\xi |^{-n++2b-|\beta |+\gamma -2\gamma_0}e^{-\sigma \rho (t,x-\xi )}.
\end{equation}

\item[(iii)]
If $|\beta |=2b,$ then
\begin{equation}
\left| D_x^\beta V_Z(t,x;\xi )\right| \le Ct^{-\alpha +\mu_1}|x-\xi |^{-n+\mu_2}e^{-\sigma \rho (t,x-\xi )};
\end{equation}
\begin{equation}
\left| D_x^\beta V_Y(t,x;\xi )\right| \le Ct^{-1+\mu_1}|x-\xi |^{-n+\mu_2}e^{-\sigma \rho (t,x-\xi )}
\end{equation}
where $\mu_1,\mu_2>0$.
\end{description}

(b) If the functions $u_0(x)$, $f(t,x)$ are bounded and globally H\"older continuous (for $f$, uniformly with respect to $t$),
and $f$ is continuous in $t$ uniformly with respect to $x\in \Rn$,
then the solution (2.25) is a uniform classical solution. All its derivatives in $x$, up to the order $2b$, are bounded and globally H\"older continuous, uniformly with respect to $t\in [0,T]$.
\end{teo}

\medskip
Note that the estimates in Theorem 2 can be written in a variety of ways. For example, in (2.26) we may write
$$
t^{-\frac{\alpha}{2b}|\beta |}=\left( t^{-\frac{\alpha}{2b}}|x-\xi |\right)^{|\beta |}|x-\xi |^{-|\beta |}
$$
and obtain, taking $0<\sigma'<\sigma$, that
$$
\left| D_x^\beta V_Z(t,x;\xi )\right| \le Ct^{-\frac{\alpha}{2b}\gamma_0}|x-\xi |^{-n-|\beta |+\gamma -\gamma_0}e^{-\sigma' \rho (t,x-\xi )}.
$$
This kind of transformation is often used in proofs of various estimates in this paper.

\bigskip
{\it 2.3. Uniqueness theorem}. Here we maintain the same assumptions as in Theorem 2.

\medskip
\begin{teo}
Let $u(t,x)$, $0\le t\le T$, $x\in \Rn$, be a uniform classical solution of the problem (1.1), (2.23) with $f(t,x)\equiv 0$, $u_0(x)\equiv 0$. Suppose that the function $u(t,x)$ and all its derivatives of orders $\le 2b$ are bounded and globally H\"older continuous. Then $u(t,x)$ equals zero identically.
\end{teo}

\medskip
The rest of the paper is devoted to the proofs of the above results. Some of them can be extended easily to more general situations. For example, the coefficients of the subordinate operator $A_1$ may depend on $t$. Some of the estimates of Theorem 1 (except the case $n+|\beta |=2b$) remain valid for systems generalizing in an obvious way the first order systems parabolic in the sense of Petrowsky. However for the whole range of the above results we need the condition (1.5).

\section{Some auxiliary results}

{\it 3.1. The Mittag-Leffler type functions of a matrix}. Let $B$ be a complex $N\times N$ matrix. The Mittag-Leffler function $E_\alpha$ of the matrix $B$ is defined by substituting $B$ into the power series (2.5):
\begin{equation}
E_\alpha (B)=\sum\limits_{k=0}^\infty \frac{B^k}{\Gamma (1+\alpha k)}, \quad 0<\alpha <1.
\end{equation}
Note that the analytic function $E_\alpha$ of a matrix does not coincide with the matrix formed by values of the function $E_\alpha$ on matrix elements.

For a class of matrices, we find a matrix analog of the asymptotic representation of the Mittag-Leffler function (see e.g. \cite{Dzh}).

\medskip
\begin{prop}
Suppose that for any $z\in \CN$,
\begin{equation}
\R \langle Bz,z\rangle \le -\delta |z|^2,\quad \delta >0.
\end{equation}
Then
\begin{equation}
E_\alpha (B)=-\frac1{\Gamma (1-\alpha )}B^{-1}+H
\end{equation}
where the matrix $H$ is such that $|H|\le C\delta^{-2}$ (the constant $C$ does not depend on $B$).
\end{prop}

\medskip
{\it Proof}. Denote by $\gamma (r,\omega )$ the contour in the complex plane oriented in the direction of the increase of $\arg \zeta$ and consisting of the following parts: the rays $\gamma_\pm =\{ \zeta \in \mathbb C:\ \arg \zeta =\pm \omega,|\zeta |\ge r\}$ and the arc $\{ \zeta \in \mathbb C:\ -\omega <\arg \zeta <\omega ,|\zeta |=r\}$. Here $r>0$, $\frac{\pi }2<\omega \le \pi$.

Let us use Hankel's integral representation
\begin{equation}
\frac1{\Gamma (s)}=\frac1{2\pi i}\int\limits_{\gamma (r,\omega )}e^\zeta \zeta^{-s}\,d\zeta ,\quad \R s>0
\end{equation}
(see e.g. \cite{Mar}). In the integral (3.4), we make the change of variables $\zeta =\eta^{1/\alpha}$. Since $r>0$ is arbitrary, we obtain the representation
\begin{equation}
\frac1{\Gamma (s)}=\frac1{2\pi i\alpha }\int\limits_{\gamma (r,\beta )}e^{\eta^{1/\alpha}} \eta^{-\frac{s}\alpha +\frac1\alpha -1}\,d\eta,
\end{equation}
for any $\beta$ with $\frac{\pi \alpha}2<\beta \le \pi \alpha$. For our purposes, we will assume that $\frac{\pi \alpha}2<\beta <\min (\frac{\pi}2,\pi \alpha )$, so that the contour $\gamma (r,\omega )$ is located in the right half-plane, and $\cos \frac{\beta}{\alpha}<0$.

Under the condition (3.2), the resolvent $(B-\lambda I)^{-1}$ exists for $\R \lambda >-\delta$, and
$$
|(B-\lambda I)^{-1}|\le \frac1{\R \lambda +\delta}
$$
(see Lemma V.6.1 in \cite{GK}). In particular, $|B^{-1}|\le \delta^{-1}$. Substituting (3.4) into (3.1) we find that
\begin{multline*}
E_\alpha (B)=\frac1{2\pi i\alpha }\sum\limits_{k=0}^\infty \left\{
\int\limits_{\gamma (r,\beta )}e^{\eta^{1/\alpha}} \eta^{-k-1}\,d\eta \right\} B^k\\
=\frac1{2\pi i\alpha }\int\limits_{\gamma (r,\beta )}e^{\eta^{1/\alpha}} \eta^{-1}\left\{ \sum\limits_{k=0}^\infty (B\eta^{-1})^k\right\} \,d\eta =-\frac1{2\pi i\alpha }\int\limits_{\gamma (r,\beta )}e^{\eta^{1/\alpha}}(B-\eta I)^{-1}\, d\eta ,
\end{multline*}
if $r>|B|$. Note that on the rays contained in  $\gamma (r,\beta )$, $\left| e^{\eta^{1/\alpha}}\right| =\exp \left( \cos \frac{\beta}{\alpha}\cdot |\eta|^{1/\alpha}\right)$ where $\cos \frac{\beta}{\alpha}<0$.

Let us fix $r=1$. The function
$$
B\mapsto \int\limits_{\gamma (1,\beta )}e^{\eta^{1/\alpha}}(B-\eta I)^{-1}\, d\eta
$$
is an analytic function of the matrix $B$ on the open set $\{ B:\
\R \langle Bz,z\rangle <-\frac{\delta}2 |z|^2,\ \forall z\in \CN \setminus \{0\} \}$ coinciding with $E_\alpha (B)$ on its intersection with the open set $\{ B:\ |B|<1\}$ (see \cite{LD} regarding analytic functions on matrices; an analytic function on $N\times N$ matrices is in fact an analytic function of $N^2$ complex variables). Therefore
\begin{equation}
E_\alpha (B)=-\frac1{2\pi i\alpha }\int\limits_{\gamma (1,\beta )}e^{\eta^{1/\alpha}}(B-\eta I)^{-1}\, d\eta
\end{equation}
for any matrix $B$ satisfying (3.2).

Next we use the identity
$$
(B-\eta I)^{-1}=B^{-1}+\eta B^{-1}(B-\eta I)^{-1}.
$$
Substituting it into (3.6) and using (3.5) we get
$$
E_\alpha (B)=-\frac1{\Gamma (1-\alpha )}B^{-1}-\frac{B^{-1}}{2\pi i\alpha }\int\limits_{\gamma (1,\beta )}e^{\eta^{1/\alpha}}\eta (B-\eta I)^{-1}\, d\eta ,
$$
which implies the required representation (3.3). $\qquad \blacksquare$

\medskip
We will need also and estimate of $E_{\alpha ,\alpha}(B)$ where the Mittag-Leffler type function $E_{\alpha ,\alpha}$ is defined by the series
\begin{equation}
E_{\alpha ,\alpha} (\zeta )=\sum\limits_{k=0}^\infty \frac{\zeta^k}{\Gamma (\alpha +\alpha k)}.
\end{equation}

\medskip
\begin{prop}
If a matrix $B$ satisfies (3.2), then
\begin{equation}
\left| E_{\alpha ,\alpha}(B)\right| \le C\delta^{-2}
\end{equation}
where $C$ does not depend on $B$.
\end{prop}

\medskip
{\it Proof}. As in the proof of Proposition 1, we obtain the representation (with the same notations)
$$
E_{\alpha ,\alpha}(B)=-\frac1{2\pi i\alpha }\int\limits_{\gamma (1,\beta )}e^{\eta^{1/\alpha}}\eta^{-1+\frac1\alpha }(B-\eta I)^{-1}\, d\eta ,
$$
from which we get that
\begin{equation}
E_{\alpha ,\alpha}(B)=-\frac1{2\pi i\alpha }B^{-1}\int\limits_{\gamma (1,\beta )}e^{\eta^{1/\alpha}}\eta^{-1+\frac1\alpha}\,d\eta
-\frac1{2\pi i\alpha }B^{-1}\int\limits_{\gamma (1,\beta )}e^{\eta^{1/\alpha}}\eta^{-1+\frac1\alpha }(B-\eta I)^{-1}\, d\eta .
\end{equation}

In the representation (3.5), we may pass to the limit, as $s\to 0$. As a result, the first integral in (3.9) equals 0. Estimating the second integral we come to (3.8). $\qquad \blacksquare$

In another result of this kind, we deal with the Mittag-Leffler type function
$$
E_{\alpha ,0} (\zeta )=\sum\limits_{k=0}^\infty \frac{\zeta^k}{\Gamma (\alpha k)}.
$$

\medskip
\begin{prop}
Under the assumptions of Proposition 1,
$$
E_{\alpha ,0}(B)=-\frac1{\Gamma (-\alpha )}B^{-1}+H,
$$
where $|H|\le C\delta^{-2}$, and the constant does not depend on $B$.
\end{prop}

\medskip
The proof is similar to that of Proposition 1.

\bigskip
{\it 3.2. Green matrices of some elliptic systems}. We will need, as technical tools, estimates of derivatives $D^\beta G(x)$ of the Green matrix of the elliptic operator $A_0(D)-I$, and also estimates of differences $D^\beta (G(x;y')-G(x;y''))$ of derivatives of the Green matrices of the operators $A_0(y',D_x)-I$ and $A_0(y'',D_x)-I$, with coefficients ``frozen'' at the points $y'$ and $y''$. Specifically, we need to consider the case where $n+|\beta |=2b$. The presence of the term $-I$ is essential -- the matrix $A_0(\eta )-I$ ($\eta \in \Rn$) has, under the assumption (2.2), all the eigenvalues with nonzero real parts. For this class of elliptic systems with constant coefficients, Eidelman \cite{E} found an integral representation of Green matrices, and the estimate
\begin{equation}
\left| D^\beta G(x)\right| \le C\left( 1+\log \frac1{|x|}\right) ,\quad |x|\le 1,
\end{equation}
where it is assumed that $n+|\beta |=2b$ (of course, other cases were considered in \cite{E} too). It also follows from the constructions in \cite{E} that $G\in C^\infty (\Rn \setminus \{0\})$.

However the estimates for the differences are not given in \cite{E}, and for completeness we have now to apply the method from \cite{E} to this situation.

Let $N(y,t,x)$ be a FSCP of the parabolic system
$$
\frac{\partial u(t,x)}{\partial t}=\left( A_0(D_x)-I\right) u(t,x).
$$
Then \cite{E}
$$
G(x;y)=\int\limits_0^\infty N(y,t,x)\,dt,
$$
so that
$$
D_x^\beta [G(x;y')-G(x;y'')]=\int\limits_0^\infty D_x^\beta [N(y',t,x)-N(y'',t,x)]\,dt.
$$

For the FSCP $N(y,t,x)$ we have the representation
$$
N(y,t,x)=(2\pi )^{-n}\int\limits_{\Rn}e^{ix\cdot \xi }e^{(A_0(y,\xi )-I)t}\,d\xi =(2\pi )^{-n}t^{-\frac{n}{2b}}e^{-t}\int\limits_{\Rn}e^{it^{-1/2b}x\cdot \xi }e^{A_0(y,\xi )}\,d\xi .
$$
Therefore, if $n+|\beta |=2b$, then
\begin{equation}
D_x^\beta [N(y',t,x)-N(y'',t,x)]=(2\pi )^{-n}t^{-1}e^{-t}\int\limits_{\Rn}\xi^\beta e^{it^{-1/2b}x\cdot \xi }\left[ e^{A_0(y',\xi )}-e^{A_0(y'',\xi )}\right]\,d\xi .
\end{equation}

Following \cite{E} we use the identity
$$
e^{A_0(y',\xi )t}-e^{A_0(y'',\xi )t}=\int\limits_0^te^{A_0(y'',\xi )(t-\tau )}[A_0(y',\xi )-A_0(y'',\xi )]e^{A_0(y',\xi )\tau }\,d\tau ,
$$
which implies the inequality
$$
\left| e^{A_0(y',\xi )}-e^{A_0(y'',\xi )}\right| \le C|y'-y''|^\gamma e^{-\delta |\xi |^{2b}}.
$$

We can consider inequalities of this kind containing, instead of $\xi \in \Rn$, a point $\xi +i\eta$, $\eta \in \Rn$. By a lemma from \cite{E} (Chapter 1), the inequality (1.5) implies the inequality
$$
\R \langle A_0(y,\xi +i\eta )z,z\rangle \le \left( -\delta_1|\xi |^{2b}+\mu_1|\eta |^{2b}\right) |f|^2
$$
($\delta_1,\mu_1>0$), thus the inequality
$$
\left| e^{A_0(y,\xi +i\eta )}\right| \le Ce^{-\delta_2|\xi |^{2b}+\mu_2|\eta |^{2b}}
$$
($\delta_2,\mu_2>0$). Repeating the above reasoning we get the inequality
\begin{equation}
\left| \xi^\beta \left[ e^{A_0(y',\xi +i\eta )}-e^{A_0(y'',\xi +i\eta )}\right] \right| \le C|y'-y''|^\gamma e^{-\delta_3|\xi |^{2b}+\mu_3|\eta |^{2b}}
\end{equation}
($\delta_3,\mu_3>0$).

Next we use another lemma from \cite{E} regarding the Fourier transform of an entire function of several complex variables satisfying an exponential estimate like (3.12). From this lemma, (3.12), and (3.11) we obtain the estimate
$$
\left| D_x^\beta [N(y',t,x)-N(y'',t,x)]\right| \le Ct^{-1}e^{-t}|y'-y''|^\gamma e^{-ct^{-\frac1{2b-1}}|x|^{\frac{2b}{2b-1}}},
$$
so that
\begin{multline*}
\left| D_x^\beta [G(x,y')-G(x,y'')]\right| \le C|y'-y''|^\gamma \int\limits_0^\infty t^{-1}e^{-t}\exp \left(-ct^{-\frac1{2b-1}}|x|^{\frac{2b}{2b-1}}\right)\,dt\\
=C|y'-y''|^\gamma \int\limits_0^\infty s^{-1}\exp \left( -s|x|^{2b}-s^{-\frac1{2b-1}}\right) \,ds=C|y'-y''|^\gamma [I_1(x)+I_2(x)]
\end{multline*}
where
$$
I_1(x)=\int\limits_0^1 s^{-1}\exp \left( -s|x|^{2b}-s^{-\frac1{2b-1}}\right) \,ds \le \int\limits_0^1 s^{-1}\exp \left( -s^{-\frac1{2b-1}}\right) \,ds\le C,
$$
\begin{multline*}
I_2(x)=\int\limits_1^\infty s^{-1}\exp \left( -s|x|^{2b}-s^{-\frac1{2b-1}}\right) \,ds\le \int\limits_1^\infty s^{-1}e^{-s|x|^{2b}}\,ds=|x|^{2b}\int\limits_1^\infty \log s\cdot e^{-s|x|^{2b}}\,ds\\
\le x|^{2b}\int\limits_0^\infty |\log s|\cdot e^{-s|x|^{2b}}\,ds =
\int\limits_0^\infty \left| \log (|x|^{-2b}\sigma )\right| e^{-\sigma }\,d\sigma \le C(1+|\log |x|).
\end{multline*}

This results in the required estimate
\begin{equation}
\left| D_x^\beta [G(x,y')-G(x,y'')]\right| \le C|y'-y''|^\gamma \left( 1+\log \frac1{|x|}\right),\quad |x|\le 1,\ n+|\beta |=2b.
\end{equation}

\bigskip
{\it 3.3. Properties of the functions $\varphi_{t,\alpha }$ and $\psi_{t,\alpha }$}. The function $\Phi_\alpha$ (see (2.4)) involved in the subordination representation (2.3) is such the $\Phi_\alpha (t)\ge 0$ for all $t>0$,
$$
\Phi_\alpha (t)\sim Ct^{-1/2}e^{-ct^{\frac1{1-\alpha }}},\quad t \to +\infty
$$
($C,c>0$); see \cite{Ba,Wr}. Changing constants if necessary, we can write the estimate
\begin{equation}
0\le \varphi_{t,\alpha }(s)\le Ct^{-\alpha }e^{-cs^{\frac1{1-\alpha }}t^{-\frac{\alpha }{1-\alpha }}},\quad s>0.
\end{equation}
Note that $\varphi_{t,\alpha }(s)\to 0$, as $t\to 0$, for each $s>0$.
It is also important that
\begin{equation}
\int\limits_0^\infty \Phi_\alpha (s)\,ds=1.
\end{equation}

In order to study the kernel $Y_\alpha$ (see (2.8)), we need the function
$$
\psi_{t,\alpha }(s)=\mathbb D_t^{1-\alpha }\varphi_{t,\alpha }(s).
$$
Making a change of variables in (2.6) we can write
$$
E_\alpha (-\zeta )=\int\limits_0^\infty \varphi_{t,\alpha }(\sigma )e^{-\zeta \sigma t^{-\alpha }}\,d\sigma ,
$$
whence
\begin{equation}
E_\alpha (-\zeta t^\alpha )=\int\limits_0^\infty \varphi_{t,\alpha }(\sigma )e^{-\zeta \sigma }\,d\sigma .
\end{equation}

The identity (3.16) is written in part as a warning of a complicated nature of the subordination identities. As $t\to +0$, the left-hand side of (3.16) tends to 1, while in the right-hand side, $\varphi_{t,\alpha }(\sigma )\to 0$ for each $\sigma >0$. Thus, it is impossible to interchange the integration and taking the limit. It will be reasonable here to use, instead of the Caputo-Dzhrbashyan derivative, the Riemann-Liouville derivative
$$
\left( D_{0+}^{1-\alpha }u\right) (t)=\frac1{\Gamma (\alpha )}\frac{d}{dt}\int\limits_0^t (t-\tau )^{-1+\alpha }u(\tau )\,d\tau
$$
coinciding with $\mathbb D^{1-\alpha }u$ wherever $u(0)=0$. In particular, if $s>0$, then
\begin{equation}
\psi_{t,\alpha }(s)=D_{0+,t}^{1-\alpha }\varphi_{t,\alpha }(s).
\end{equation}

By (2.4), we have the Wright function representation
$$
\varphi_{t,\alpha }(s)=t^{-\alpha }
{}_0\Psi_1\Bigl[ \begin{matrix}
- \\ (1-\alpha ,-\alpha )\end{matrix}\Bigl| -st^{-\alpha } \Bigr] .
$$
Using the contour integral representation of the Wright function (see the equality (12.41) in \cite{B}) we find that
\begin{equation}
\varphi_{t,\alpha }(s)=\frac1{2\pi i}\int\limits_{\gamma -i\infty }^{\gamma +i\infty }\frac{\Gamma (\lambda )}{\Gamma (1-\alpha +\alpha \lambda )}s^{-\lambda }t^{\alpha \lambda -\alpha }\,d\lambda ,
\end{equation}
$\gamma >0$, $\gamma \ne \frac{\alpha +\nu }\alpha$, $\nu =0,1,2,\ldots$.

It is known (see, for example, \cite{KST}) that $D_{0+}^{1-\alpha }$ transforms $t^{\alpha \lambda -\alpha }$ into $\dfrac{\Gamma (1-\alpha +\alpha \lambda )}{\Gamma (\alpha \lambda )}t^{\alpha \lambda -1}$. Now we get from (3.17) and (3.18) that
$$
\psi_{t,\alpha }(s)=\frac1{2\pi i}\int\limits_{\gamma -i\infty }^{\gamma +i\infty }\frac{\Gamma (\lambda )}{\Gamma (\alpha \lambda )}s^{-\lambda }t^{\alpha \lambda -1}\,d\lambda
$$
(the legitimacy of this transformation follows from the asymptotics of Gamma function). Using again the equivalence of the series and integral representations of the Wright functions ((12.41) in \cite{B}) we prove the representaton (2.9):
\begin{equation}
\psi_{t,\alpha }(s)=t^{-1}
{}_0\Psi_1\Bigl[ \begin{matrix}
- \\ (0,-\alpha )\end{matrix}\Bigl| -st^{-\alpha } \Bigr] ,
\end{equation}
so that
\begin{equation}
\psi_{t,\alpha }(s)=t^{-1}\sum\limits_{k=1}^\infty \frac{(-st^{-\alpha } )^k}{k!\Gamma (-\alpha k)}.
\end{equation}

Using the asymptotics of the Wright function found in \cite{Wr} (see also \cite{B}, Theorem 25) we get the estimate
\begin{equation}
\left| \psi_{t,\alpha }(s)\right| \le Ct^{-1}\exp \left\{ -c(st^{-\alpha })^{\frac{1}{1-\alpha }}\right\},
\end{equation}
with $c,C>0$.

Following (2.6) and (3.16), let us write a connection of the function $\psi_{t,\alpha }$ with the Mittag-Leffler type functions. Apply the operator $D_{0+,t}^{1-\alpha }$ to both sides of the equality (3.16). It is known (see (2.1.54) in \cite{KST}) that in the left-hand side we obtain $t^{\alpha -1}E_{\alpha ,\alpha }(-\zeta t^\alpha )$. Thus we come to the identity
\begin{equation}
E_{\alpha ,\alpha }(-\zeta t^\alpha )=t^{-\alpha }
\int\limits_0^\infty {}_0\Psi_1\Bigl[ \begin{matrix}
- \\ (0,-\alpha )\end{matrix}\Bigl| -st^{-\alpha } \Bigr]
e^{-\zeta s}\,ds,\quad \zeta >0.
\end{equation}

Another transformation kernel of the above kind is
\begin{equation}
\nu_{t,\alpha }(s)=\frac{\partial}{\partial t} \varphi_{t,\alpha }(s).
\end{equation}
Repeating the above reasoning we find that
\begin{equation}
\nu_{t,\alpha }(s)=t^{-\alpha -1}
{}_0\Psi_1\Bigl[ \begin{matrix}
- \\ (-\alpha ,-\alpha )\end{matrix}\Bigl| -st^{-\alpha } \Bigr]
\end{equation}
where the Wright function has the form
$$
{}_0\Psi_1\Bigl[ \begin{matrix}
- \\ (-\alpha ,-\alpha )\end{matrix}\Bigl| z \Bigr] =\sum\limits_{k=0}^\infty \frac{(-z)^k}{k!\Gamma (-\alpha -\alpha k)}.
$$

Using the asymptotics from \cite{B,Wr} we obtain, for $s>0$, $t>0$, the estimate
\begin{equation}
\left| \nu_{t,\alpha }(s)\right| \le Ct^{-\alpha -1}\exp \left\{ -c(st^{-\alpha })^{\frac{1}{1-\alpha }}\right\},\quad c>0.
\end{equation}

It is known (\cite{KST}, formula (1.10.2)) that
$$
\frac{\partial}{\partial t}E_\alpha (-\zeta t^\alpha )=t^{-1}E_{\alpha ,0}(-\zeta t^\alpha )
$$
where $E_{\alpha ,0}(\zeta )=\sum\limits_{k=0}^\infty \dfrac{\zeta^k}{\Gamma (\alpha k)}.
$
Differentiating both sides of the identity (3.16) we find, after elementary transformations, that
\begin{equation}
E_{\alpha ,0}(-\zeta )=
\int\limits_0^\infty {}_0\Psi_1\Bigl[ \begin{matrix}
- \\ (-\alpha ,-\alpha )\end{matrix}\Bigl| -s\Bigr]
e^{-\zeta s}\,ds,\quad \zeta >0.
\end{equation}

\section{Systems with constant coefficients}

{\it 4.1. Proof of Theorem 1}. Let us use the identity (2.7) which implies a representation of spatial derivatives,
\begin{equation}
D_x^\beta Z_\alpha (t,x)=\int\limits_0^\infty \varphi_{t,\alpha }(s)D_x^\beta Z(s,x)\,ds,\quad x\ne 0.
\end{equation}
Below we use both the integral representation for $Z$,
\begin{equation}
Z(s,x)=(2\pi )^{-n}\int\limits_{\Rn}e^{ix\cdot \xi }e^{sA_0(\xi )}\,d\xi ,
\end{equation}
and its estimates
\begin{equation}
\left| D_x^\beta Z(t,x)\right| \le Ct^{-\frac{n+|\beta |}{2b}}\exp \left\{ -c|x|^{\frac{2b}{2b-1}}t^{-\frac1{2b-1}}\right\}
\end{equation}
($c>0$) valid for all $x\in \Rn$, $t>0$ (see \cite{E}).

Note that, by its construction, the function $Z_\alpha (t,x)$ satisfies the system (2.1) in a weak sense (a strong solution is obtained after a convolution in spatial variables with a test function; see Section 2.1). It will be seen from the investigation of the spatial derivatives below that $A_0(D_x)Z_\alpha$ exists in the classical sense for $t>0$ and $x\ne 0$. Therefore $\left( \D Z_\alpha \right) (t,x)$ coincides almost everywhere with $A_0(D_x)Z_\alpha (t,x)$ ($t>0,x\ne 0$) and does not require a separate investigation.

It follows from the inequalities (3.14) and (4.3) that $Z_\alpha (t,x)$ is infinitely differentiable in $x\ne 0$, and
$$
\left| D_x^\beta Z_\alpha (t,x)\right| \le Ct^{-\alpha } \int\limits_0^\infty s^{-\frac{n+|\beta |}{2b}}\exp \left\{ -c|x|^{\frac{2b}{2b-1}}s^{-\frac1{2b-1}}\right\} \exp \left\{ -cs^{\frac1{1-\alpha }}t^{-\frac{\alpha }{1-\alpha }}\right\} \,ds
$$
(for $t>0$, the second exponential factor on the right guarantees the convergence of the integral at infinity; for $x\ne 0$, the first exponential factor gives the convergence at the origin).
After the change of variables $s^{-\frac1{2b-1}}=\sigma$, we find that
\begin{equation}
\left| D_x^\beta Z_\alpha (t,x)\right| \le Ct^{-\alpha } \int\limits_0^\infty \sigma^{\frac{2b-1}{2b}(n+|\beta |)-2b}\exp \left\{ -c\sigma |x|^{\frac{2b}{2b-1}}\right\}\exp \left\{ -c\sigma^{-\frac{2b-1}{1-\alpha }}t^{-\frac{\alpha }{1-\alpha }}\right\} \,d\sigma .
\end{equation}

In order to obtain an estimate for $R\ge 1$, we use the asymptotics of the integral
\begin{equation}
\Omega (\zeta )=\int\limits_0^\infty e^{-\zeta t}e^{-dt^{-\varkappa }}t^\lambda \,dt,\quad \zeta \to \infty ,
\end{equation}
where $d>0$, $\varkappa >0$, $\lambda \in \mathbb R$, found in \cite{Ri} (formula (12.80)). Namely,
\begin{equation}
\Omega (\zeta )\sim a_0(d\varkappa \zeta^{-1})^{\frac{\lambda +1}{\varkappa +1}}\exp \left[ -(1+\frac1\varkappa )\rho \right] \rho^{-1/2}
\end{equation}
where $\rho =(d\varkappa \zeta^\varkappa )^{\frac1{1+\varkappa }}$, $a_0=2\left( \frac2{1+\varkappa }\right)^{1/2}\Gamma (\frac12)$.

Making in (4.4) the change of variables $\sigma =t^{-\frac{\alpha }{2b-1}}\eta$ we get, after easy calculations, the inequality
$$
\left| D_x^\beta Z_\alpha (t,x)\right| \le Ct^{-\alpha \frac{n+|\beta |}{2b} } \int\limits_0^\infty \eta^{\frac{2b-1}{2b}(n+|\beta |)-2b} \exp\left( -c\eta^{-\frac{2b-1}{1-\alpha}}\right) \exp \left( -cR^\frac1{2b-1}\eta \right)\,d\eta
$$
where the integral has the form of (4.5) with $\zeta =R^\frac1{2b-1}$, $d=c$, $\varkappa =\dfrac{2b-1}{1-\alpha}$, $\lambda =\dfrac{2b-1}{2b}(n+|\beta |)-2b$. Using (4.6) and ignoring powers of $R$ (changing, if necessary, the constant in the exponential factor) we obtain the inequality (2.10).

Suppose that $R\le 1$, $n+|\beta |>2b$. We use again the inequality (4.4), but make the change of variables $\sigma =\tau |x|^{-\frac{2b}{2b-1}}$ and replace the exponential factor containing $|x|$ by 1. We find that
$$
\left| D_x^\beta Z_\alpha (t,x)\right| \le Ct^{-\alpha } |x|^{-n-|\beta |+2b} \int\limits_0^\infty \tau^{\frac{2b-1}{2b}(n+|\beta |)-2b}e^{-c\tau }\,d\tau
$$
where $\dfrac{2b-1}{2b}(n+|\beta |)-2b>-1$. This implies (2.14).

The inequalities (2.11) and (2.15) are proved similarly, on the basis of the estimate (3.21).

Let $R\le 1$, $n+|\beta |<2b$. Performing a change of variables we can rewrite (4.2) in the form
$$
Z(s,x)=(2\pi )^{-n}s^{-\frac{n}{2b}}\int\limits_{\Rn}e^{is^{-1/2b}x\cdot \xi }e^{A_0(\xi )}\,d\xi ,
$$
and by virtue of (4.1),
$$
D_x^\beta Z_\alpha (t,x)=(2\pi )^{-n}t^{-\alpha }\int\limits_0^\infty \Phi_\alpha (st^{-\alpha })s^{-\frac{n+|\beta |}{2b}}\,ds\int\limits_{\Rn}\xi^\beta e^{is^{-1/2b}x\cdot \xi }e^{A_0(\xi )}\,d\xi ,
$$
that is, after the change $s=\sigma t^\alpha$,
\begin{equation}
D_x^\beta Z_\alpha (t,x)=(2\pi )^{-n}t^{-\alpha \frac{n+|\beta |}{2b}}\int\limits_0^\infty \Phi_\alpha (\sigma )\sigma^{-\frac{n+|\beta |}{2b}}\,d\sigma \int\limits_{\Rn}\xi^\beta e^{it^{-\alpha /2b}\sigma^{-1/2b}x\cdot \xi }e^{A_0(\xi )}\,d\xi .
\end{equation}

It follows from (2.2) that
\begin{equation}
\left| e^{A_0(\xi )}\right| \le e^{-\delta |\xi |^{2b}}
\end{equation}
(see Sect. I.4.4 in \cite{Kr}). Using the fact that $\Phi_\alpha$ decays rapidly at infinity, and $\dfrac{n+|\beta |}{2b}<1$, we come to the inequality (2.12). The proof of (2.13), based on (2.8), (3.19), and (3.21), is similar.

Finally, consider the most complicated case where $R\le 1$, $n+|\beta |=2b$. The representation (4.7) takes the form
\begin{equation}
D_x^\beta Z_\alpha (t,x)=(2\pi )^{-n}t^{-\alpha}\int\limits_0^\infty \Phi_\alpha (\sigma )\sigma^{-1}\,d\sigma \int\limits_{\Rn}\xi^\beta e^{it^{-\alpha /2b}\sigma^{-1/2b}x\cdot \xi }e^{A_0(\xi )}\,d\xi
\end{equation}
or, after the change $\xi =\sigma^{1/2b}\eta$,
\begin{equation}
D_x^\beta Z_\alpha (t,x)=(2\pi )^{-n}t^{-\alpha}\int\limits_0^\infty \Phi_\alpha (\sigma )\,d\sigma \int\limits_{\Rn}\eta^\beta e^{it^{-\alpha /2b}x\cdot \eta }e^{\sigma A_0(\eta )}\,d\eta .
\end{equation}

The identity (2.6) remains valid when the matrix $-A_0(\eta )$ is substituted for $\zeta$. Indeed, we may rewrite (2.6) in the form
$$
E_\alpha (-\zeta )=\sum\limits_{k=0}^\infty \frac{(-1)^k\zeta^k}{k!}\int\limits_0^\infty \Phi_\alpha (t)t^k\,dt
$$
(the convergence of the series follows from the asymptotics of $\Phi_\alpha$ and Stirling's formula). Recall that an entire function of a matrix is defined by substituting the matrix into the power series expansion. Thus we set $\zeta =-A_0(\eta )$, and then gather the power series into the exponential (keeping in mind the inequality (4.8)). Therefore
\begin{equation}
E_\alpha (A_0(\eta ))=\int\limits_0^\infty \Phi_\alpha (\sigma )e^{\sigma A_0(\eta )}\,d\sigma.
\end{equation}

If we substitute (4.11) into (4.10) and change the order of integration, we obtain the representation
\begin{equation}
D_x^\beta Z_\alpha (t,x)=(2\pi )^{-n}t^{-\alpha}\int\limits_{\Rn}\eta^\beta e^{it^{-\alpha /2b}x\cdot \eta }E_\alpha (A_0(\eta ))\,d\eta .
\end{equation}
However this change of the order of integration requires a justification.

It follows from (2.2) and Proposition 1 that
$$
\left| E_\alpha (A_0(\eta ))\right| \le C|\eta |^{-2b},\quad |\eta |\ge 1,
$$
so that $|\eta |^{|\beta |}\left| E_\alpha (A_0(\eta ))\right| \le C|\eta |^{-n}$, $|\eta |\ge 1$. Denote
$$
X_\varepsilon (t,x,\beta )=(2\pi )^{-n}t^{-\alpha}\int\limits_\varepsilon^\infty \Phi_\alpha (\sigma )\sigma^{-1}\,d\sigma \int\limits_{\Rn}\xi^\beta e^{it^{-\alpha /2b}\sigma^{-1/2b}x\cdot \xi }e^{A_0(\xi )}\,d\xi ,\quad \varepsilon >0.
$$

Here we make the change of variables $\xi =\sigma^{1/2b}\eta$, and change the order of integration (we have moved away from the singularity!). Thus
$$
X_\varepsilon (t,x,\beta )=(2\pi )^{-n}t^{-\alpha}\int\limits_{\Rn}\eta^\beta e^{it^{-\alpha /2b}x\cdot \eta }\,d\eta \int\limits_\varepsilon^\infty \Phi_\alpha (\sigma )e^{\sigma A_0(\eta )}d\sigma .
$$
Let us show that
\begin{equation}
\eta^\beta \int\limits_\varepsilon^\infty \Phi_\alpha (\sigma )e^{\sigma A_0(\eta )}d\sigma \xrightarrow{\varepsilon \to 0}\eta^\beta \int\limits_0^\infty \Phi_\alpha (\sigma )e^{\sigma A_0(\eta )}d\sigma =E_\alpha (A_0(\eta ))\eta^\beta ,
\end{equation}
as functions of $\eta$, in the topology of $L_2(\Rn )$.

Indeed, using the estimate (4.8) we see that
\begin{multline*}
\int\limits_{\Rn}|\eta |^{2|\beta |}\left| \int\limits_0^\varepsilon \Phi_\alpha (\sigma )e^{\sigma A_0(\eta )}d\sigma \right|^2\,d\eta \le C\int\limits_{\Rn}|\eta |^{2|\beta |}\left[ \int\limits_0^\varepsilon e^{-\delta \sigma |\eta |^{2b}}d\sigma \right]^2\,d\eta \\
\le \int\limits_{\Rn}|\eta |^{-2n}\left( 1-e^{-\delta \varepsilon |\eta |^{2b}}\right)^2\,d\eta =C(I_1+I_2)
\end{multline*}
where
$$
I_1=\int\limits_{|\eta |\le 1}|\eta |^{-2n}\left( 1-e^{-\delta \varepsilon |\eta |^{2b}}\right)^2\,d\eta ,
$$
$$
\int\limits_{|\eta |>1}|\eta |^{-2n}\left( 1-e^{-\delta \varepsilon |\eta |^{2b}}\right)^2\,d\eta .
$$

Using the inequality $1-e^{-x}\le x$, $x\ge 0$, we find that
$$
I_1\le C\varepsilon^2\int\limits_{|\eta |\le 1}|\eta |^{-2n+4b}\,d\eta =C\varepsilon^2\int\limits_{|\eta |\le 1}|\eta |^{2|\beta |}\,d\eta\to 0,
$$
as $\varepsilon \to 0$. By the dominated convergence theorem, we obtain also that $I_2\to 0$, and we have proved (4.13).

Now, by the properties of the Fourier transform, for any fixed $t>0$ and almost all $x\in \Rn$,
$$
\lim\limits_{\varepsilon \to 0}X_\varepsilon (t,x,\beta )=(2\pi )^{-n}t^{-\alpha}\int\limits_{\Rn}\eta^\beta e^{it^{-\alpha /2b}x\cdot \eta }E_\alpha (A_0(\eta ))\,d\eta .
$$
On the other hand, $X_\varepsilon (t,x,\beta )\to D_x^\beta Z_\alpha (t,x)$, if $x\ne 0$, by (4.9). Thus we have proved the equality (4.12) for almost all $x\ne 0$.

Denote temporarily the right-hand side of (4.12) by $X(t,x,\beta )$. Next we prove that $X(t,x,\beta )$ is continuous in $x\ne 0$. This will establish the equality (4.12) for all $x\ne 0$; simultaneously we will get the required estimate.

By Proposition 1, $X(t,x,\beta )=X_1+X_2+X_3$ where
$$
X_1=(2\pi )^{-n}t^{-\alpha}\int\limits_{|\eta |\le 1}\eta^\beta e^{it^{-\alpha /2b}x\cdot \eta }E_\alpha (A_0(\eta ))\,d\eta ,
$$
$$
X_2=-\frac{t^{-\alpha }}{(2\pi )^n\Gamma (1-\alpha )}\int\limits_{|\eta |>1}\eta^\beta e^{it^{-\alpha /2b}x\cdot \eta }[A_0(\eta )]^{-1}\,d\eta ,
$$
$$
X_3=(2\pi )^{-n}t^{-\alpha}\int\limits_{|\eta |>1}\eta^\beta e^{it^{-\alpha /2b}x\cdot \eta }H(\eta )\,d\eta ,
$$
$|H(\eta )|\le C|\eta |^{-4b}$. Since $|\beta |-4b=-n-2b$, $X_3$ is continuous in $x$, and $|X_3|\le Ct^{-\alpha }$. We see also that $X_1$ is continuous in $x$, and $|X_1|\le Ct^{-\alpha }$.

Let us write $X_2=X_{21}+X_{22}+X_{23}$ where
$$
X_{21}=-\frac{t^{-\alpha }}{(2\pi )^n\Gamma (1-\alpha )}\int\limits_{|\eta |>1}\eta^\beta e^{it^{-\alpha /2b}x\cdot \eta }\left\{ [A_0(\eta )]^{-1}-[A_0(\eta )-I]^{-1}\right\} \,d\eta ,
$$
$$
X_{22}=-\frac{t^{-\alpha }}{(2\pi )^n\Gamma (1-\alpha )}\int\limits_{\Rn }\eta^\beta e^{it^{-\alpha /2b}x\cdot \eta }
[A_0(\eta )-I]^{-1}\,d\eta ,
$$
$$
X_{23}=-\frac{t^{-\alpha }}{(2\pi )^n\Gamma (1-\alpha )}\int\limits_{|\eta |\le 1}\eta^\beta e^{it^{-\alpha /2b}x\cdot \eta }[A_0(\eta )-I]^{-1}\,d\eta .
$$

We have
\begin{equation}
\left| [A_0(\eta )-I]^{-1}\right| \le \frac1{1+\delta |\eta |^{2b}}
\end{equation}
(see Lemma V.6.1 in \cite{GK}). By the resolvent identity
$$
[A_0(\eta )]^{-1}-[A_0(\eta )-I]^{-1}=[A_0(\eta )]^{-1}[A_0(\eta )-I]^{-1},
$$
we get the estimate
\begin{equation}
\left| [A_0(\eta )]^{-1}-[A_0(\eta )-I]^{-1}\right| \le \frac1{\delta|\eta |^{2b}(1+\delta |\eta |^{2b})}.
\end{equation}
It follows from (4.14) and (4.15) that $X_{21}$ and $X_{23}$ are continuous in $x$, $|X_{21}|\le Ct^{-\alpha }$, $|X_{23}|\le Ct^{-\alpha }$.

As for $X_{22}$, we note that
$$
X_{22}=-\frac{t^{-\alpha }}{\Gamma (1-\alpha )}D_y^\beta G(y)\left|
{}_{y=t^{-\alpha /2b}x}\right.
$$
where $G$ is the Green matrix of the elliptic operator $A_0(D)-I$. This means that $X_{22}$ is continuous in $x\ne 0$, thus $X$ has the same property and coincides with $D^\beta Z_\alpha$ for all $x\ne 0$. Now the required inequality (2.17) is a consequence of (3.10).

For $n=1$, we will refine this estimate. In this case, $A_0(\eta )=a_0\eta^{2b}$ where
$$
\R \langle a_0z,z\rangle \le -\delta |z|^2\quad \text{for all $z\in \CN$}.
$$
We have
$$
D_x^\beta Z_\alpha (t,x)=\frac1{2\pi }t^{-\alpha }\int\limits_{-\infty}^\infty \eta^\beta e^{it^{-\alpha /2b}x\eta }E_\alpha (a_0\eta^{2b})\,d\eta
$$
where $1+\beta =2b$. In particular, the natural number $\beta$ is odd, so that
$$
D_x^\beta Z_\alpha (t,x)=\frac1{\pi }t^{-\alpha }\int\limits_0^\infty \eta^\beta \sin (t^{-\alpha /2b}x\eta )E_\alpha (a_0\eta^{2b})\,d\eta .
$$

By Proposition 1,
$$
E_\alpha (a_0\eta^{2b})=-\frac{a_0^{-1}}{\Gamma (1-\alpha )}\eta^{-2b}+O(\eta^{-4b}),\quad \eta \to \infty .
$$
The contribution of the remainder term in the estimate of $D_x^\beta Z_\alpha$ is clearly $O(t^{-\alpha})$. Therefore we have to consider the function
\begin{equation}
F(y)=\int\limits_0^\infty \varphi (\eta )\sin (y\eta )\,d\eta ,\quad 0<y\le 1,
\end{equation}
where $\varphi$ is continuous on $[0,\infty )$, $\varphi (\eta )\sim \eta^{-1}$, $\eta \to \infty$.

It was shown in \cite{SS} that the integral in (4.16) exists as an improper one, and
$$
F(y)\sim \const \cdot y^{-1}\varphi (y^{-1}),\quad y\to 0
$$
(asymptotics of this kind is proved in \cite{SS} for much more general situations), that is $F$ is bounded near the origin. This implies the inequality (2.16).

For the function $Y_\alpha$ with $n+|\beta |=2b$, we use, in a similar way, the identity (3.22) with a matrix argument $\zeta$, which results in the representation
$$
D_x^\beta Y_\alpha (t,x)=(2\pi )^{-n}t^{-1}\int\limits_{\Rn}\xi^\beta e^{it^{-\alpha /2b}x\cdot \xi }E_{\alpha ,\alpha }(A_0(\xi ))\,d\xi .
$$
Using Proposition 2 we obtain the desired estimate (2.18).

The estimates for the first time derivative $\dfrac{\partial}{\partial t}Z_\alpha (t,x)$ are obtained just as those for $Z_\alpha$ itself. We use the representations (3.23), (3.24), and (3.26), as well as the estimate (3.25) and the matrix asymptotics given by Proposition 3. With this input, the proofs of (2.19)-(2.22) are similar to the ones given above. $\qquad \blacksquare$

\bigskip
{\it 4.2. Unified estimates}. In Theorem 1, the estimates are given separately for large and small values of $R$. In order to justify the iteration procedures of the Levi method, we need unified estimates valid for all values of the variables.

\medskip
\begin{prop}
If $n+|\beta |<2b$, then
\begin{equation}
\left| D_x^\beta Z_\alpha (t,x)\right| \le Ct^{-\alpha \frac{n+|\beta |}{2b}}e^{-c\rho (t,x)},\quad c>0;
\end{equation}
\begin{equation}
\left| D_x^\beta Y_\alpha (t,x)\right| \le Ct^{-1+\alpha -\alpha \frac{n+|\beta |}{2b}}e^{-c\rho (t,x)}.
\end{equation}

If $n+|\beta |>2b$, then
\begin{equation}
\left| D_x^\beta Z_\alpha (t,x)\right| \le Ct^{-\alpha}|x|^{-n+2b-|\beta|}e^{-c\rho (t,x)}.
\end{equation}

If $n+|\beta |=2b$, then
\begin{equation}
\left| D_x^\beta Z_\alpha (t,x)\right| \le Ct^{-\alpha}\left[ \left| \log \left( t^{-\alpha }|x|^{2b}\right) \right| +1\right] e^{-c\rho (t,x)}.
\end{equation}

If $n+|\beta |\ge 2b$, then
\begin{equation}
\left| D_x^\beta Y_\alpha (t,x)\right| \le Ct^{-1}|x|^{-n+2b-|\beta|} e^{-c\rho (t,x)}.
\end{equation}
The constants can depend only on the parameters listed in the formulation of Theorem 1.
\end{prop}

\medskip
{\it Proof}. The estimate (4.17) coincides with (2.10), if $R\ge 1$, being obviously equivalent to (2.12), if $R\le 1$, $n+|\beta |<2b$.

Let us consider the case where $n+|\beta |>2b$. It is clear that (4.19) is equivalent to (2.14), if $R\le 1$. If $R\ge 1$, we rewrite the right-hand side of (2.10) as follows. Let $\sigma =c'+c''$ ($c',c''>0$). Then
\begin{multline*}
t^{-\alpha \frac{n+|\beta |}{2b}}e^{-c(t^{-\alpha /2b}|x|)^{\frac{2b}{2b-\alpha }}}\\
=t^{-\alpha}|x|^{-n+2b-|\beta|}\left[ \left( t^{-\alpha /2b}|x|\right)^{n+|\beta|-2b}e^{-c'(t^{-\alpha /2b}|x|)^{\frac{2b}{2b-\alpha }}}\right] e^{-c''(t^{-\alpha /2b}|x|)^{\frac{2b}{2b-\alpha }}}\\
\le Ct^{-\alpha}|x|^{-n+2b-|\beta|}e^{-c\rho (t,x)}
\end{multline*}
where $c=c''$, and we have proved (4.19).

The proofs of (4.18), (4.20), and (4.21) are similar. $\qquad \blacksquare$

\section{Parametrix}

The parametrix kernels $Z_\alpha^{(0)}(t,x;y)$ and $Y_\alpha^{(0)}(t,x;y)$ defined in Section 2.2 satisfy all the estimates of Theorem 1 and Proposition 4, with all the constants independent of $y$.

We need also estimates of the differences $Z_\alpha^{(0)}(t,x;y')-Z_\alpha^{(0)}(t,x;y'')$ and $Y_\alpha^{(0)}(t,x;y')-Y_\alpha^{(0)}(t,x;y'')$. These estimates are identical to those for $Z_\alpha^{(0)}$ and $Y_\alpha^{(0)}$ themselves, with an additional factor $|y'-y''|^\gamma$. The proofs are the same as in Theorem 1 and Proposition 4, with the following additional tools: the difference estimates for classical parabolic systems \cite{E}; the difference estimate (3.13) for the Green matrices of elliptic systems; the estimate for
$$
\left| e^{A_0(y',\eta )}- e^{A_0(y'',\eta )}\right|
$$
given in Chapter 1 of \cite{E}. We omit further details since they just repeat the above material.

As in Section 2.1, we have the integral identities
\begin{equation}
\int\limits_{\Rn} Z_\alpha^{(0)}(t,x;y)\,dx=1,\quad \int\limits_{\Rn} Y_\alpha^{(0)}(t,x;y)\,dx=\frac{t^{\alpha -1}}{\Gamma (\alpha)}.
\end{equation}
It follows from the difference estimates and the first identity in (5.1) that
\begin{equation}
\left| \int\limits_{\Rn}\frac{\partial }{\partial t} Z_\alpha^{(0)}(t,x-\xi ;\xi )\,d\xi \right| \le Ct^{-1+\frac{\alpha \gamma}{2b}}.
\end{equation}

\section{The Levi method. Proof of Theorem 2}

Given the estimates of Theorem 1 and Proposition 4, the proof of Theorem 2 is carried out just as its counterpart for fractional diffusion equations \cite{EK,EIK}. The integral inequalities needed for the proof are given in sufficient generality in \cite{EIK}. Therefore we drop the detailed calculations and give only the scheme and the main estimates.

We look for the functions $Z_\alpha^{(1)}(t,x;\xi )$, $Y_\alpha^{(1)}(t,x;\xi )$ appearing in Theorem 2 assuming the following integral representations:
\begin{equation}
Z_\alpha^{(1)}(t,x;\xi )=Z_\alpha^{(0)}(t,x-\xi ;\xi )+\int\limits_0^td\lambda
\int\limits_{\mathbb R^n}Y_\alpha^{(0)}(t-\lambda ,x-y;y)Q(\lambda ,y;\xi )\,dy;
\end{equation}
\begin{equation}
Y_\alpha^{(1)}(t,x;\xi )=Y_\alpha^{(0)}(t,x-\xi ;\xi )+\int\limits_0^td\lambda
\int\limits_{\mathbb R^n}Y_\alpha^{(0)}(t-\lambda ,x-y;y)\Phi (\lambda ,y;\xi )\,dy.
\end{equation}

For the functions $Q,\Phi$ we assume the integral equations
\begin{equation}
Q(t,x;\xi )=M(t,x;\xi )+\int\limits_0^td\lambda
\int\limits_{\mathbb R^n}K(t-\lambda,x;y)Q(\lambda ,y;\xi
)\,dy,
\end{equation}
\begin{equation}
\Phi (t,x;\xi )=K(t,x;\xi )+\int\limits_0^td\lambda
\int\limits_{\mathbb R^n}K(t-\lambda,x;y)\Phi (\lambda ,y;\xi )\,dy,
\end{equation}
where
$$
M(t,x;\xi )=[A(x,D_x)-A_0(\xi ,D_x)]Z_\alpha^{(0)}(t,x-\xi ;\xi
),
$$
$$
K(t,x;\xi )=[A(x,D_x)-A_0(\xi ,D_x)]Y_\alpha^{(0)}(t,x-\xi ;\xi
)
$$

Using the estimates from Proposition 4 we find that
\begin{equation}
|M(t,x;\xi )|\le Ct^{-\alpha }|x-\xi |^{-n+\gamma }e^{-c\rho (t,x-\xi )},
\end{equation}
\begin{equation}
|K(t,x;\xi )|\le Ct^{-1+(\gamma-\eta )\frac{\alpha }{2b}}|x-\xi |^{-n+\eta }e^{-c\rho (t,x-\xi )},
\end{equation}
$c>0$, $0<\eta <\gamma$.

The increments of $M$ and $K$ are estimated as follows. Let $\Delta_xM(t,x;\xi )=M(t,x;\xi )-M(t,x';\xi )$,
$\Delta_xK(t,x;\xi )=K(t,x;\xi )-K(t,x';\xi )$. Denote by $x''$ one of the points $x,x'$, for which $|x''-\xi |=\min \{ |x-\xi |,|x'-\xi |\}$. Then
\begin{equation}
|\Delta_xM(t,x;\xi )|\le Ct^{-\alpha }|x-x'|^{\gamma -\varepsilon
}|x''-\xi |^{-n+\varepsilon }\exp \{-\sigma \rho (t,x''-\xi )\},
\end{equation}
\begin{equation}
|\Delta_xK(t,x;\xi )|\le Ct^{-1}|x-x'|^{\gamma -\varepsilon
}|x''-\xi |^{-n+\varepsilon }\exp \{-\sigma \rho (t,x''-\xi
)\},
\end{equation}
$\varepsilon >0$.

Using the estimates (6.5)-(6.8) we prove the convergence of iterations and obtain estimates for the solutions of the integral equations (6.3) and (6.4). We use, as a tool, Lemma 1.14 from \cite{EIK} (where the conditions $n=2$ and $\nu_0<1$ are in fact unnecessary). The resulting estimates are as follows:
\begin{equation}
|Q(t,x;\xi )|\le Ct^{-\alpha }|x-\xi |^{-n+\gamma }\exp
\{-\sigma \rho (t,x-\xi )\},
\end{equation}
\begin{equation}
|\Phi (t,x;\xi )|\le Ct^{-1}|x-\xi |^{-n+\gamma }\exp
\{-\sigma \rho (t,x-\xi )\},
\end{equation}
\begin{equation}
|\Delta_xQ(t,x;\xi )|\le Ct^{-\alpha }|x-x'|^{\gamma -\varepsilon
}|x''-\xi |^{-n+\varepsilon }\exp \{-\sigma \rho (t,x''-\xi )\},
\end{equation}
\begin{equation}
|\Delta_x\Phi (t,x;\xi )|\le Ct^{-1}|x-x'|^{\gamma -\varepsilon
}|x''-\xi |^{-n+\varepsilon }\exp \{-\sigma \rho (t,x''-\xi )\},
\end{equation}

Now the representation of the Green matrix stated in Theorem 2 follows from its construction (6.1)-(6.2) while the estimates (2.26)-(2.31) are obtained from (6.9)-(6.12) and Lemmas 1.12, 1.13 from \cite{EIK}. The above estimates, together with the inequality (5.2), make it possible also to repeat, without significant changes, the whole reasoning from \cite{EK} or \cite{EIK} regarding the heat potential and the initial condition.

It also follows from (6.11), (6.12), and the difference estimates of Section 5, that the differences $D_x^\beta \left[ V_Z(t,x';\xi )-V_Z(t,x'';\xi )\right]$, $D_x^\beta \left[ V_Y(t,x';\xi )-V_Y(t,x'';\xi )\right]$, $|\beta |\le 2b$, satisfy the estimates similar to those for $D_x^\beta V_Z$, $D_x^\beta V_Y$, with the additional factor $|x'-x''|^{\gamma -\varepsilon }$.

Let us find out when the solution $u(t,x)$ of the form (2.25) is a uniform classical solution with the global H\"older properties. We consider the more complicated second term
$$
w(t,x)=\int\limits_0^td\lambda \int\limits_{\mathbb R^n} Y_\alpha^{(1)}(t-\lambda,x;y)f(\lambda ,y)\,dy=w_1(t,x)+w_2(t,x)
$$
where
$$
w_1(t,x)=\int\limits_0^td\lambda \int\limits_{\mathbb R^n} Y_\alpha^{(0)}(t-\lambda,x-y;y)f(\lambda ,y)\,dy,
$$
$$
w_2(t,x)=\int\limits_0^td\lambda \int\limits_{\mathbb R^n} V_Y(t-\lambda,x;y)f(\lambda ,y)\,dy
$$
(the first term in (2.25) can be considered similarly). The remainder kernel $V_Y$ is less singular than $Y_\alpha^{(0)}$, so that the uniform convergence of derivatives of $w_2$ is verified in a straightforward way. The same can be said about lower order derivatives of $w_1$.

For the leading derivatives, we have, just as in \cite{EK} or \cite{EIK}, the expressions
\begin{multline}
D_x^\beta w_1(t,x)=\int\limits_0^td\lambda \int\limits_{\mathbb R^n} D_x^\beta Y_\alpha^{(0)}(t-\lambda,x-y;y)[f(\lambda ,y)-f(\lambda ,x)]\,dy\\
+\int\limits_0^tf(\lambda ,x)\,d\lambda \int\limits_{\mathbb R^n} D_x^\beta Y_\alpha^{(0)}(t-\lambda,x-y;y)\,dy,\quad |\beta |=2b,
\end{multline}
\begin{multline}
\left( \D w_1\right) (t,x)=f(t,x)+\int\limits_0^td\lambda
\int\limits_{\mathbb R^n}
\frac{\partial Z_\alpha^{(0)}(t-\lambda ,x-y;y)}{\partial t}[f(\lambda ,y)-f(\lambda
,x)]\,dy\\
+\int\limits_0^tf(\lambda ,x)\,d\lambda \int\limits_{\mathbb R^n}
\frac{\partial Z_\alpha^{(0)}(t-\lambda ,x-y;y)}{\partial t}\,dy.
\end{multline}

The global H\"older property of the derivatives (6.13) follows from the difference estimates of $D_x^\beta Y_\alpha^{(0)}$ and our assumptions regarding the function $f$.

The representation (6.14) is obtained as follows (see \cite{EK} or \cite{EIK}). First of all, if $v(t,x)=\left( I_{0+}^{1-\alpha }w_1\right) (t,x)$, then $\D w_1=\dfrac{\partial v}{\partial t}$, and
$$
v(t,x)=\int\limits_0^td\lambda \int\limits_{\mathbb R^n}
Z_\alpha^{(0)}(t-\lambda ,x-y;y)f(\lambda ,y)\,dy.
$$

For a small positive number $h$, set
$$
v_h(t,x)=\int\limits_0^{t-h}d\lambda \int\limits_{\mathbb R^n}
Z_\alpha^{(0)}(t-\lambda ,x-y;y)f(\lambda ,y)\,dy.
$$
Then $\dfrac{\partial v_h}{\partial t}=v_h^{(1)}+v_h^{(2)}$ where
$$
v_h^{(1)}(t,x)=\int\limits_{\mathbb R^n}
Z_\alpha^{(0)}(h,x-y;y)f(t-h,y)\,dy,
$$
$$
v_h^{(2)}(t,x)=\int\limits_0^{t-h}d\lambda \int\limits_{\mathbb R^n}
\frac{\partial Z_\alpha^{(0)}(t-\lambda ,x-y;y)}{\partial t}f(\lambda ,y)\,dy.
$$
We have
\begin{multline*}
v_h^{(1)}(t,x)=\int\limits_{\mathbb R^n}
\left[ Z_\alpha^{(0)}(h,x-y;y)-Z_\alpha^{(0)}(h,x-y;x)\right] f(t-h,y)\,dy\\
+\int\limits_{\mathbb R^n}Z_\alpha^{(0)}(h,x-y;x)[f(t-h,y)-f(t-h,x)]\,dy+f(t-h,x).
\end{multline*}

It follows from the estimates of the parametrix kernel and its differences, and from the global H\"older property of $f$ that both integrals in the last formula tend to zero, as $h\to 0$, uniformly with respect to $x\in \Rn$, $t\in [0,T]$.

Similarly, we prove the convergence of $v_h^{(2)}(t,x)$ to the sum of the two integrals appearing in (6.14), uniformly with respect to $x,t$. This proves the uniform property of our solution. $\qquad \blacksquare$

\section{Proof of Theorem 3}

First we consider the model system (2.1) with constant coefficients. A uniform classical solution of the system (2.1) can be interpreted as a classical solution of the operator-differential equation
\begin{equation}
\left( \D w\right) (t)=B_0w(t)
\end{equation}
in the Banach space $C_b(\Rn )^N$ of bounded continuous vector-functions with the supremum norm. Here $B_0$ is the closed operator on $C_b(\Rn )^N$ defined by $A_0(D)$ with the domain
$$
\left\{ v\in C_b(\Rn )^N:\ A_0(D)v\in C_b(\Rn )^N\right\}
$$
($A_0(D)v$ is understood in the sense of tempered distributions). Let $q$ be such a natural number that $q>(N+\frac{n}2)\cdot 2b$. By Theorem 4.1 of the paper \cite{HHN}, under the condition (2.2) (in fact, the Petrowsky parabolicity condition would suffice), we have
\begin{equation}
\left\| (\lambda I-B_0)^{-1}(I-\Delta)^{-q/2}\right\| \le p(|\lambda |),\quad \R \lambda >0,
\end{equation}
where $p$ is a certain polynomial.

Note that the operator $(I-\Delta)^{-q/2}$ is bounded on $C_b(\Rn )^N$; that follows from the integrability of its integral kernel \cite{St}. Under the assumption (7.2), the equation (7.1) has only a trivial solution $w\in C_b(\Rn )^N$ with $w(0)=0$. That is proved exactly as the uniqueness theorem from \cite{K1} where it was assumed that $\limsup\limits_{\lambda \to \infty}\lambda^{-1/\alpha}\log \left\| (\lambda I-B_0)^{-1}\right\| =0$. One should just repeat the whole reasoning from \cite{K1} for the function $(I-\Delta)^{-q/2}w$, instead of $w$, and notice that the operator $(I-\Delta)^{-q/2}$ is injective.

Thus, we have proved Theorem 3 for the system (2.1). Turning to the general case, we rewrite the system (1.1) with $f=0$ in the form
\begin{equation}
\left( \D u\right) (t,x)-A_0(y,D_x)u(t,x)=[A_0(x,D_x)-A_0(y,D_x)]u(t,x)+A_1(x,D_x)u(t,x).
\end{equation}
Here $y\in \Rn$ is an arbitrary fixed point. As before, we assume that $u(0,x)=0$.

By Theorem 2, we can write down an integral representation of a uniform classical solution of the Cauchy problem for the equation (7.3). By the above uniqueness result for model systems with constant coefficients, we obtain the equality
\begin{multline}
u(t,x)=\int\limits_0^td\tau \int\limits_{\mathbb R^n} Y_\alpha^{(0)}(t-\tau,x-\xi ;y)[A_0(\xi ,D_\xi )-A_0(y,D_\xi )]u(t,\xi )\,d\xi \\
+\int\limits_0^td\tau \int\limits_{\mathbb R^n} Y_\alpha^{(0)}(t-\tau,x-\xi ;y)A_1(\xi ,D_\xi )u(t,\xi )\,d\xi .
\end{multline}

Denote
$$
v(t)=\sum\limits_{|m|\le 2b}\sup\limits_{x\in \Rn }\left| D_x^m u(t,x)\right| .
$$
Differentiating both sides of (7.4) in $x$ and taking into account the boundedness of the derivatives of $u$ we come to the inequality
$$
\left| D_x^m u(t,x)\right| \le C\int\limits_0^tv(\tau )\,d\tau \int\limits_{\mathbb R^n}\left| D_x^mY_\alpha^{(0)}(t-\tau,x-\xi ;y)\right| \,d\xi ,\quad |m|<2b.
$$

Using the estimates (4.18) and (4.21) we find that
\begin{equation}
\left| D_x^m u(t,x)\right| \le C\int\limits_0^tv(\tau )(t-\tau )^{-1+(2b-|m|)\frac{\alpha }{2b}}\,d\tau ,\quad |m|<2b.
\end{equation}

For $|m|=2b$, the derivatives of a heat potential are regularized by subtraction, that is
\begin{multline*}
D_x^m u(t,x)=\int\limits_0^td\tau \int\limits_{\mathbb R^n} D_x^m Y_\alpha^{(0)}(t-\tau,x-\xi ;y)\left\{ [A_0(\xi ,D_\xi )-A_0(y,D_\xi )]u(t,\xi )\right. \\
\left. -[A_0(x,D_x)-A_0(y,D_x)]u(t,x)\right\} \,d\xi \\
+\int\limits_0^td\tau \int\limits_{\mathbb R^n} D_x^m Y_\alpha^{(0)}(t-\tau,x-\xi ;y)[A_1(\xi ,D_\xi )u(t,\xi )-A_1(x,D_x)u(t,x)]\,d\xi .
\end{multline*}
Until now, $y\in \Rn$ was an arbitrary parameter. Set $y=x$. Then one of the terms in the first integral disappears. For the remaining terms, we use the inequalities
$$
|[A_0(\xi ,D_\xi )-A_0(x,D_\xi )]u(t,\xi )| \le C|x-\xi |^\gamma v(t),
$$
$$
|A_1(\xi ,D_\xi )u(t,\xi )-A_1(x,D_x)u(t,x)| \le C|x-\xi |\sup_{\substack{x\in \Rn \\ |\beta |\le \deg A_1+1}}\left| D_x^\beta u(t,x)\right| \le C|x-\xi |v(t)
$$
(since $\deg A_1+1\le 2b$).

It follows from (4.21) that
\begin{multline*}
\left| D_x^m u(t,x)\right| \le C\int\limits_0^tv(\tau )(t-\tau )^{-1}\,d\tau \\
\times \int\limits_{\mathbb R^n}(|x-\xi |^\gamma +|x-\xi |)|x-\xi |^{-n}\exp \left\{ -c\left[ (t-\tau )^{-\alpha}|x-\xi |^{2b}\right]^{\frac1{2b-\alpha }}\right\} \,d\xi \\
\le C\int\limits_0^tv(\tau )(t-\tau )^{-1+\frac{\alpha \gamma }{2b}}\,d\tau ,
\end{multline*}
$|m|=2b$. Adding this inequality to (7.5) we find that, for any $t\in [0,T]$,
\begin{equation}
v(t)\le C\int\limits_0^tv(\tau )(t-\tau )^{-1+\varepsilon}\,d\tau ,\quad \varepsilon >0.
\end{equation}
By a kind of the Bellman-Gronwall inequality proved by Henry (\cite{H}, Lemma 7.1.1), it follows from (7.6) that $v(t)\equiv 0$, hence $u(t,x)\equiv 0$. $\qquad \blacksquare$

\end{document}